\documentclass{article}    
\usepackage{amstext}
\def\qed{$\rlap{$\sqcap$}\sqcup$}
\usepackage{latexsym}

\begin{document}           

{\ }\\
\begin{center}
{\huge {\bf The $h$-vector of a relatively\\compressed level algebra}} \\ [.250in]
{\large FABRIZIO ZANELLO\\
Department of Mathematics, Royal Institute of Technology, 100 44 Stockholm, Sweden\\
E-mail: zanello@math.kth.se}
\end{center}

{\ }\\
\\
ABSTRACT. The purpose of this note is to supply an upper and a lower bound (which are in general sharp) for the $h$-vector of a level algebra which is relatively compressed with respect to any arbitrary level algebra $A$.\\
The useful concept of relatively compressed algebra was recently introduced in $[MMN]$ by Migliore {\it et al.} (whose investigations mainly focused on the particular case of $A$ a complete intersection). The key idea of this note is the simple observation that the level algebras which are relatively compressed with respect to $A$ coincide (after an obvious isomorphism) with the generic level quotients of suitable truncations of $A$. Therefore, we are able to apply to relatively compressed algebras the main result of our recent work $[Za3]$.\\
\\
\section{Introduction}
\indent
{\large

The idea of a {\it compressed algebra} has been around in commutative algebra for some time (it first appeared, implicitly, in Emsalem-Iarrobino's 1978 work $[EI]$), and identifies those (standard graded artinian) algebras having the (entry by entry) maximum $h$-vector among the $h$-vectors of all the algebras having given codimension and socle-vector. Compressed algebras and their $h$-vectors were extensively studied in the eighties by Iarrobino ($[Ia]$) and Fr\"oberg-Laksov ($[FL]$), and then by this author in full generality ($[Za1]$ and $[Za2]$).\\\indent
Recently, Migliore, Mir\'o-Roig and Nagel ($[MMN]$) defined the useful concept of {\it relatively compressed algebra} (which extends a notion first introduced by Cho and Iarrobino; see $[CI]$, Remark 2.9): given an artinian algebra $A=R/I$, an algebra $B=R/J$ (where $R$ is a polynomial ring and $I$ and $J$ are homogeneous ideals of $R$) is called {\it relatively compressed with respect to $A$}, if $B$ is a quotient algebra of $A$ (i.e., $J\supseteq I$) and the $h$-vector of $B$ is the largest among the $h$-vectors of all the quotient algebras of $A$ having the same socle-vector as $B$.\\\indent
This concept, in other words, defines the \lq \lq smallest" ideals $J$ containing a fixed ideal $I$, given the socle-vector of $R/J$ as an invariant. Migliore {\it et al.} mainly focused on studying the $h$-vectors (and the minimal free resolutions) of those level algebras which are relatively compressed with respect to $R$ modulo a complete intersection (oftentimes generic). Those $h$-vectors, in fact, supply interesting information on the size of the ideals containing a given regular sequence.\\\indent
In this note, by identifying relatively compressed algebras with generic quotients of suitable truncations of a level algebra, we can apply the recent work contained in $[Za3]$, where we studied generic level quotients, to extend some of the results of $[MMN]$. Namely, we supply (in general sharp) upper and lower bounds for the $h$-vector of a level algebra which is relatively compressed with respect to any given arbitrary level algebra.\\
\\\indent
Let us now fix the setting we will be working in throughout this note. We consider standard graded artinian algebras $A=R/I$, where $R=k[x_1,...,x_r]$, $I$ is a homogeneous ideal of $R$, $k$ is a field of characteristic zero and the $x_i$'s all have degree 1.\\\indent
The {\it $h$-vector} of $A$ is $h(A)=h=(h_0,h_1,...,h_e)$, where $h_i=\dim_k A_i$ and $e$ is the last index such that $\dim_k A_e>0$. Since we may suppose that $I$ does not contain non-zero forms of degree 1, $r=h_1$ is defined to be the {\it codimension} of $A$.\\\indent 
The {\it socle} of $A$ is the annihilator of the maximal homogeneous ideal $\overline{m}=(\overline{x_1},...,\overline{x_r})\subseteq A$, namely $soc(A)=\lbrace a\in A {\ } \mid {\ } a\overline{m}=0\rbrace $. Since $soc(A)$ is a homogeneous ideal, we define the {\it socle-vector} of $A$ as $s(A)=s=(s_0,s_1,...,s_e)$, where $s_i=\dim_k soc(A)_i$. Note that $h_0=1$, $s_0=0$ and $s_e=h_e>0$. The integer $e$ is called the {\it socle degree} of $A$ (or of $h$). The {\it type} of the socle-vector $s$ (or of the algebra $A$) is type$(s)=\sum_{i=0}^es_i$.\\\indent 
If $s=(0,0,...,0,s_e=t)$, we say that the algebra $A$ is {\it level} (of type $t$). With a slight abuse of notation, we will refer to an $h$-vector as level if it is the $h$-vector of a level algebra.\\\indent
Let us now recall the main facts of the theory of {\it inverse systems} that we will be using in this paper. For a complete introduction, we refer the reader to $[Ge]$ and $[IK]$.\\\indent 
Let $S=k[y_1,...,y_r]$, and consider $S$ as a graded $R$-module where the action of $x_i$ on $S$ is partial differentiation with respect to $y_i$.\\\indent 
There is a one-to-one correspondence between artinian algebras $R/I$ and finitely generated $R$-submodules $M$ of $S$, where $I=Ann(M)$ is the annihilator of $M$ in $R$ and, conversely, $M=I^{-1}$ is the $R$-submodule of $S$ which is annihilated by $I$ (cf. $[Ge]$, Remark 1), p. 17).\\\indent 
If $R/I$ has socle-vector $s$, then $M$ is minimally generated by $s_i$ elements of degree $i$, for $i=1,2,...,e$, and the $h$-vector of $R/I$ is given by the number of linearly independent derivatives in each degree obtained by differentiating the generators of $M$ (cf. $[Ge]$, Remark 2), p. 17).\\\indent 
In particular, level algebras of type $t$ and socle degree $e$ correspond to $R$-submodules of $S$ minimally generated by $t$ elements of degree $e$.\\
\\
\section{The main result}
\indent

Let $A=R/I$ be any level algebra of socle degree $e$, having $h$-vector $h=(1,h_1,...,h_e)$. Let $B$ be a level algebra of socle degree $d\leq e$ and type $c\leq h_d$ which is relatively compressed with respect to $A$. Denote the $h$-vector of $B$ by $H=(1,H_1,...,H_d)$. Then:\\
\\\indent
{\bf Theorem 1.} {\it With the above notations, for $i=1,2,...,d$, we have:
$${1\over h_d^2-1}((h_d-c)h_{d-i}+(ch_d-1)h_i)\leq H_i\leq \min \lbrace h_i,ch_{d-i}\rbrace .$$}
\\\indent
{\bf Proof.} Notice that a quotient of socle degree $d$ of $A$ can be naturally seen as a quotient of socle degree $d$ of $A/A_{d+1}\cong R/(I,R_{d+1})$, which is the level truncation of $A$ having $h$-vector $h=(1,h_1,...,h_d)$.\\\indent
Also note the important fact that the level algebras $B$ above, relatively compressed with respect to $A$, are exactly (isomorphic to) the {\it generic} level quotients of type $c$ and socle degree $d$ of the algebra $R/(I,R_{d+1})$ (that is, the quotients parameterized by the points of a non-empty open subset of $({\bf P}^{h_d-1}(k))^c$, given a set of generators $G_1$, ..., $G_{h_d}$ of the inverse system module corresponding to the level algebra $R/(I,R_{d+1})$; see also $[Za3]$).\\\indent
Hence, the first inequality of the theorem is immediately obtained by $[Za3]$, Theorem 2.9 (Main Theorem).\\\indent
Let us now show the second inequality: $H_i\leq h_i$ is obvious. It remains to prove that $H_i\leq ch_{d-i}$. Again, consider $B$ as a quotient algebra of $R/(I,R_{d+1})$. The inverse system module corresponding to $B$ is generated by $c$ linearly independent forms of degree $d$, and $H_i$ is given by the number of $(d-i)$-th linearly independent partial derivatives of those $c$ forms. By the symmetry of {\it Gorenstein} $h$-vectors (i.e., level $h$-vectors of type 1), each of the above $c$ forms has at most $h_{d-i}$ $(d-i)$-th linearly independent partial derivatives, whence we have $H_i\leq ch_{d-i}$. This completes the proof of the theorem.{\ }{\ }\qed \\
\\\indent 
{\bf Remark 2.} The second inequality of Theorem 1 extends $[MMN]$, Lemma 2.13 to any level algebra $A$ (not necessarily a complete intersection).\\
\\\indent
{\bf Proposition 3.} {\it Let $A$ be any level algebra having $h$-vector $h=(1,h_1,...,h_e)$. Then, for all positive integers $d$ and $c$ such that $d\leq e$ and $c\leq h_d$, there exists a level algebra having socle degree $d$ and type $c$ which is relatively compressed with respect to $A$.}\\
\\\indent
{\bf Proof.} As we have underlined in the proof of the theorem, the level algebras of socle degree $d$ which are relatively compressed with respect to a given level algebra $A=R/I$ \lq \lq are" exactly the generic level quotients of $R/(I,R_{d+1})$ of socle degree $d$. The result follows from the fact that those generic quotients always exist (see the observation at the beginning of Section 2 of $[Za3]$).{\ }{\ }\qed \\
\\\indent
{\bf Example 4.} Let us consider a level algebra $A$ of socle degree 8 having $h$-vector $$h=(1,4,9,13,13,13,9,6,4).$$\indent
(The fact that $h$ is actually a level $h$-vector can be shown, for instance, by using sums of powers of generic linear forms as suggested in $[Ia]$, Theorem 4.8 A.) By Theorem 1, the $h$-vector $H$ of the level algebras of type 3 and socle degree 6 which are relatively compressed with respect to $A$ must satisfy the following (entry by entry) inequalities: $$(1,3,4,6,5,5,3)\leq H\leq (1,4,9,13,13,12,3).$$
\\\indent
{\bf Remark 5.} In general, the two bounds of Theorem 1 are sharp (and they are so simultaneously for all $i=1,2,...,d$).\\\indent
Indeed, the fact that there exist algebras $A$ and $B$ such that the lower bound of the theorem is sharp was implicitly shown in $[Za3]$, Example 2.10. We briefly recap that example here (without proofs) for completeness.\\\indent
Let us consider the level algebra $A$ (of codimension $r=(t+1)p$, type $t$ and socle degree $e$) associated to the inverse system module $M=<F_1,...,F_t>$, where $F_j=y_{jp+1}y_1^{e-1}+y_{jp+2}y_2^{e-1}+...+y_{(j+1)p}y_p^{e-1}$. One can easily see that the $h$-vector of $A$ is $h=(1,(t+1)p,(t+1)p,...,(t+1)p,t)$. Furthermore, the level quotients of $A$ of type $c\leq t$ and socle degree $e$ all have the same $h$-vector, $(1,(c+1)p,(c+1)p,...,(c+1)p,c)$. But this is exactly the lower bound of the theorem, since we have the identity $(c+1)p={1\over t^2-1}((t-c)(t+1)p+(ct-1)(t+1)p)$. Thus, in this example the lower bound of Theorem 1 is sharp.\\\indent
The existence of algebras $A$ and $B$ such that the $h$-vector of $B$ is the upper bound of Theorem 1 is guaranteed by an argument using compressed algebras. For instance, let us consider two generic forms, $F$ and $G$, of degree 7 inside $S=k[y_1,y_2,y_3]$. Then the $h$-vector of the (compressed) level algebra $A=R/Ann(F,G)$ is $h=(1,3,6,10,15,12,6,2)$ (e.g., see $[Ia]$). If we consider the level quotient $B$ of $A$ of type 2 and socle degree 6 whose associated inverse system module is generated by two derivatives of $F$ and $G$ (e.g., $F_{y_1}$ and $G_{y_2}$), then we clearly have that the $h$-vector of $B$ is $H=(1,3,6,10,12,6,2)$, and this is exactly the upper bound supplied by Theorem 1.\\
\\\indent
{\bf Remark 6.} A very interesting fact regarding relatively compressed algebras is that their $h$-vector depends heavily on the original algebra $A$, not only on its $h$-vector.\\\indent
In order to see this, we will exhibit two level algebras, $A_1$ and $A_2$, such that $h(A_1)=h(A_2)=h=(1,h_1,...,h_e)$, and two positive integers, $d\leq e$ and $c\leq h_d$, such that there exist level algebras $B_1$ and $B_2$, both having socle degree $d$ and type $c$, which are relatively compressed with respect to $A_1$ and $A_2$ (respectively), and whose $h$-vectors are different.\\\indent
Let $h=(1,(t+1)p,(t+1)p,...,(t+1)p,t)$, as in Remark 5, and consider, as $A_1$, the algebra $A$ of that remark. Let $d=e$ and $c=1$. As we saw above, the $h$-vector of any Gorenstein quotient $B_1$ of $A_1$ of socle degree $e$ is $H^{(1)}=(1,2p,2p,...,2p,1)$.\\\indent
Therefore, it suffices to find a level algebra $A_2$ with $h$-vector $h$ having a Gorenstein quotient of socle degree $e$ with an $h$-vector $H^{(2)}$ larger than $H^{(1)}$. Let us consider an inverse system submodule $M\subset S=k[y_1,y_2,...,y_{(t+1)p}]$ generated by $t$ forms of degree $e$, one (say $F$) being the sum of powers of $(t+1)p-(t-1)$ generic linear forms, and the remaining $t-1$ being the $e$-th powers of one generic linear form each. If we let $A_2$ be the level algebra associated to $M$, then it is easy to see (again by using $[Ia]$, Theorem 4.8 A on sums of powers of linear forms) that the $h$-vector of $A_2$ is $h$, and that the $h$-vector of the Gorenstein quotient $B_2$ of $A_2$ (of socle degree $e$) whose inverse system cyclic module is generated by $F$ is $H^{(2)}=(1,(t+1)p-(t-1),(t+1)p-(t-1),...,(t+1)p-(t-1),1)$.\\\indent
Thus, for $t>1$ and $p>1$, we have $H^{(2)}>H^{(1)}$, as we desired.\\
\\\indent
{\bf Remark 7.} It would be interesting to extend the above results on level algebras to algebras $A$ and $B$ with arbitrary socle-vectors. However, this is beyond the scope of this brief note, since such a generalization would be strongly related to an analogous generalization of our recent paper $[Za3]$, that, to date, has not yet been obtained.\\
\\\indent
{\bf Acknowledgements.} We warmly thank (our former Ph.D. Thesis advisor) Dr. A.V. Geramita and Dr. A. Iarrobino for their comments and suggestions which improved the exposition of this article.\\
\\
\\
{\bf \huge References}\\
$[CI]$ {\ } Y.H. Cho and A. Iarrobino: {\it Hilbert Functions and Level Algebras}, J. of Algebra 241 (2001), 745-758.\\
$[EI]$ {\ } J. Emsalem and A. Iarrobino: {\it Some zero-dimensional generic singularities; finite algebras having small tangent space}, Compositio Math. 36 (1978), 145-188.\\
$[FL]$ {\ } R. Fr\"oberg and D. Laksov: {\it Compressed Algebras}, Conference on Complete Intersections in Acireale, Lecture Notes in Mathematics, No. 1092 (1984), 121-151, Springer-Verlag.\\
$[Ge]$ {\ } A.V. Geramita: {\it Inverse Systems of Fat Points: Waring's Problem, Secant Varieties and Veronese Varieties and Parametric Spaces of Gorenstein Ideals}, Queen's Papers in Pure and Applied Mathematics, No. 102, The Curves Seminar at Queen's (1996), Vol. X, 3-114.\\
$[Ia]$ {\ } A. Iarrobino: {\it Compressed Algebras: Artin algebras having given socle degrees and maximal length}, Trans. Amer. Math. Soc. 285 (1984), 337-378.\\
$[IK]$ {\ } A. Iarrobino and V. Kanev: {\it Power Sums, Gorenstein Algebras, and Determinantal Loci}, Springer Lecture Notes in Mathematics (1999), No. 1721, Springer, Heidelberg.\\
$[MMN]$ {\ } J. Migliore, R. Mir\'o-Roig and U. Nagel: {\it Minimal resolution of relatively compressed level algebras}, J. of Algebra 284 (2005), No. 1, 337-370.\\
$[Za1]$ {\ } F. Zanello: {\it Extending the idea of compressed algebra to arbitrary socle-vectors}, J. of Algebra 270 (2003), No. 1, 181-198.\\
$[Za2]$ {\ } F. Zanello: {\it Extending the idea of compressed algebra to arbitrary socle-vectors, II: cases of non-existence}, J. of Algebra 275 (2004), No. 2, 730-748.\\
$[Za3]$ {\ } F. Zanello: {\it Partial derivatives of a generic subspace of a vector space of forms: quotients of level algebras of arbitrary type}, Trans. Amer. Math. Soc., to appear.

}

\end{document}